\documentclass[12pt]{article}
\usepackage{amssymb,amsmath,theorem}
 \let\kappa=\varkappa                  %  греческие буквы
\let\emptyset=\varnothing                                  %  пустое множество
\newcommand\qed{\ifhmode\unskip\nobreak\fi\quad            %  квадрат в конце
   \ifmmode\square\else\hbox{$\square$}\fi}                %   доказательств
\newcommand\pin{\kern.0833em}                              %  половина тонкого пробела
\theorempreskipamount=\medskipamount \theorempostskipamount=\medskipamount

 \newtheorem{theorem}{Theorem}%[section]
 \newtheorem{lemma}[theorem]{Lemma}
 \newtheorem{corollary}[theorem]{Corollary}
{\theorembodyfont{\normalfont}
 \newtheorem{remark}[theorem]{Remark}
 
\newtheorem{definition}{Definition}
}
\begin{document}

\begin{center} \bf \Large
Sup\kern0.03em -sums principles for \emph{F}-divergence,\\ Kullback---Leibler divergence, and\\ new definition for \emph{t}-entropy

\bigskip\medskip\normalsize\rm

V.\,I.\ BAKHTIN

\smallskip

{\it John Paul II Catholic University of Lublin, Poland\ \ {\rm\&}\ \

Belarusian State University, Belarus

$($e-mail: bakhtin@tut.by\/$)$}

\bigskip

A.\,V.\ LEBEDEV

\smallskip

{\it University of Bialystok, Poland\ \ {\rm\&}\ \

Belarusian State University, Belarus

$($e-mail: lebedev@bsu.by\/$)$}
\end{center}

\vspace{-18pt}

\renewcommand\abstractname{}
\begin{abstract} \noindent
The article presents new   ${\sup}$-sums principles for integral $F$-divergence for arbitrary convex function $F$ and arbitrary (not necessarily positive and absolutely continuous) measures. As applications of these results we derive the corresponding   ${\sup}$-sums principle for  Kullback---Leibler divergence and work out new `integral' definition for $t$-entropy explicitly establishing its relation to Kullback---Leibler divergence.
\end{abstract}

\bigbreak\bigskip

\quad\parbox{13.9cm} {\textbf{Keywords:} {\itshape $F$-divergence, Kullback---Leibler divergence, sup-sums principle, partition of unity, $t$-entropy}

\medbreak

}

\section*{Introduction}

The notion of $F$-divergence was introduced and originally studied in analysis of probability distributions by Csisz\'ar (1963), Morimoto (1963) and Ali, Silvey (1966) \cite{Csis, Mor,AlSil}. It is defined in the following way. Let $P$ and $Q$ be two probability distributions over a space $\Omega $ such that $P$ is absolutely continuous with respect to $Q$. Then, for a convex function $F\!: \mathbb R_+ \to \mathbb R$ such that $F(1) = 0$, the $F$-divergence $D_F (P\Vert Q)$ of $P$ from $Q$ is defined as
\begin{equation*}
 D_F (P\Vert Q) := \int_\Omega F\bigg(\frac{d\,P}{d\,Q}\bigg)\, dQ,
\end{equation*}
where $dP/dQ$ is the Radon---Nikodym derivative of $P$ with respect to $Q$.

Since its introduction $F$-divergence has been intensively exploited and analysed due to the fact that by taking appropriate functions $F$ one arrives here at numerous important divergences such as Kullback---Leibler divergence, Hellinger distance, Pearson $\chi^2$-di\-ver\-gence etc.

In the article we derive a number of new  ${\sup}$-sums principles for integral $F$-divergence type objects that concern not only probability distributions but also real-valued measures and are associated with general convex functions $F$ that can take infinite values (Theorems~\ref{..theorem_1a},  \ref{..theorem_3} and  \ref{..theorem_1}). These  ${\sup}$-sums principles express integrals in question by means of integral sums type objects associated with partitions of unity. In particular, they give us an explicit relation between $F$-divergences for continuous and discrete measures.
One can also note here a new observation of the arising in these principles additional parameters $F'(\pm \infty)$.

On the base of   ${\sup}$-sums principles obtained we derive the corresponding
${\sup}$-sums principle for Kullback---Leibler divergence (Theorem~\ref{t.KL}) leading also naturally to its new definition for measures that are not probability ones.

As one more substantial application of integral  ${\sup}$-sums principles deduced we worked out  a principally new definition for $t$-entropy. The $T$-entropy plays a fundamental role in the spectral analysis of operators associated with dynamical systems (cf. Theorem~\ref{1..1}) and so also is a key ingredient in `entropy statistic theorem' --- the statement that plays in the spectral theory of weighted shift and transfer operators the role analogous to Shannon---McMillan---Breiman theorem in information theory \cite{McMill,AlCov} and its important corollary known as `asymptotic equipartition property' \cite[p.~135]{Bill}. Up to now the definition of $t$-entropy has been formulated in a rather sophisticated manner in terms of actions of transfer operators on continuous partitions of unity (detailed discussion see in Subsection~\ref{s-n-t-entr}). In Theorem~\ref{t-t-entr-KL} we give a fundamentally new `integral' definition of $t$-entropy explicitly establishing its relation to Kullback---Leibler divergence.

\section{Sup\kern0.03em -sums \emph{F}-divergence} \label{s-0}

Consider an arbitrary convex function $F\!:\mathbb R\to (-\infty, +\infty]$. Let
\begin{equation} \label{,,71}
 F'(+\infty) :=\lim_{t\to +\infty}\frac{F(t)}{t}, \qquad
 F'(-\infty) :=\lim_{t\to -\infty} \frac{F(t)}{t}.
\end{equation}
Obviously, both limits do exist, and the value of $F'(+\infty)$ may be finite or equal to $+\infty$ while $F'(-\infty)$ may be finite or equal to $-\infty$.

Now we adopt the following agreement. The product $0 F(x/0)$ for $x\ne 0$ will be defined as limit $\lim_{t\to +0} tF(x/t)$, and for $x=0$ it will be assumed to be zero. In other words,
\begin{equation} \label{,,72}
 0F\bigg(\frac{x}{0}\bigg) \,=\,
  \begin{cases}
    \pin xF'(+\infty), & \text{if}\ \, x>0, \\[2pt]
    \pin xF'(-\infty), & \text{if}\ \, x<0, \\[2pt]
    \pin 0,            & \text{if}\ \, x=0.
  \end{cases}
\end{equation}

Let a finite nonnegative measure $\mu$ and a finite real-valued measure $\nu$ be defined on a measurable space $(X,\mathfrak A)$. Then for measurable functions $g$ on $(X,\mathfrak A)$ we will exploit the notation
\begin{equation*}
 \mu[g] := \int_X g\,d\mu, \qquad \nu[g] := \int_X g\,d\nu
\end{equation*}

 \medskip\noindent
(under assumption that the integrals do converge).

By a \emph{measurable partition of unity} we will understand a finite set $G = \{g_1,\dots,g_k\}$ of nonnegative measurable functions on $(X,\mathfrak A)$ such that $\sum_i g_i \equiv 1$.

%\begin{definition} \label{d-,,rho}
For any convex function $F\!:\mathbb R \to (-\infty,+\infty]$ set
\begin{equation} \label{,,rho}
 \rho^{}_F(\mu,\nu) :=\, \sup_G\sum_{g\in G} \mu[g]\pin F\bigg(\frac{\nu[g]}{\mu[g]}\bigg),
\end{equation}
where supremum is taken over the set of all measurable partitions of unity\/ $G$ and we assume that if\/ $\mu[g]=0$ then the corresponding summand in the right-hand part is defined according to convention \eqref{,,72}.

Bearing in mind the initial flavor of $F$-divergence it is natural to call $\rho^{}_F(\mu,\nu)$ as \emph{sup-sums $ F$-divergence}. Its relation to the usual (integral) $ F$-divergence will be uncovered in the next section.
%\end{definition}

The principal role in the definition \eqref{,,rho} is played by the function $sF(x/s)$. The next two lemmas describe its important technical properties that also will be exploited in the proofs of further results.

\begin{lemma} \label{..72}
For any convex function\/ $F$ and all\/ $s,t\ge 0$ and\/ $x,y\in\mathbb R,$
\begin{equation} \label{,,75}
 (s+t)F\bigg(\frac{x+y}{s+t}\bigg) \pin\le\pin sF\bigg(\frac{x}{s}\bigg) +\pin tF\bigg(\frac{y}{t}\bigg).
\end{equation}
\end{lemma}

Each convex function $F$ on the real axis is superlinear, i.\,e.,
\begin{equation} \label{,,73}
 F(t) \ge At +B
\end{equation}
for some constants $A,B\in\mathbb R$ and all $t\in\mathbb R$.

\begin{lemma} \label{..71}
If a convex function\/ $F$ satisfies condition\/ \eqref{,,73}, then for all\/ $s\ge 0$ and\/ $x\in\mathbb R,$
\begin{equation} \label{,,74}
 sF\bigg(\frac{x}{s}\bigg) \ge Ax +Bs.
\end{equation}
\end{lemma}

Now we proceed to description of the principal technical properties of $\rho^{}_F(\mu,\nu)$.

\begin{lemma} \label{..lemma_rho}
Expression\/ \eqref{,,rho} is well-defined and preserves its value would we use countable partitions of unity\/ $G$ in it instead of finite ones.
\end{lemma}

\begin{theorem} \label{..73}
The function\/ $\rho^{}_F(\mu,\nu)$ is subadditive with respect to the pair\/ $(\mu,\nu)$. That is, for any finite nonnegative measures\/ $\mu_1,\,\mu_2$ and any finite real-valued measures\/ $\nu_1,\,\nu_2,$
\begin{equation} \label{,,76}
 \rho^{}_F(\mu_1+\mu_2,\nu_1+\nu_2) \pin\le\pin \rho^{}_F(\mu_1,\nu_1) +\rho^{}_F(\mu_2,\nu_2).
\end{equation}
\end{theorem}

For any measure $\nu$ and bounded measurable function $f$ on a measurable space $(X,\mathfrak A)$ one can define a real-valued measure $f\nu$ by the rule
\begin{equation*}
 f\nu[g] := \nu[f g] =\int_X gf\,d\nu, \qquad g\in L^1(X,\nu).
\end{equation*}

\begin{theorem} \label{..additivity}
Let\/ $\mu,\,\nu$ be finite measures, where\/ $\mu$ is nonnegative and\/ $\nu$ is real-valued, and\/ $f_1,\pin f_2$ be nonnegative bounded measurable functions on\/ $(X,\mathfrak A)$. Then
\begin{equation} \label{,,additivity}
 \rho^{}_F\big((f_1+f_2)\mu,(f_1+f_2)\nu\big) =\rho^{}_F(f_1\mu,f_1\nu) +\rho^{}_F(f_2\mu,f_2\nu).
\end{equation}

 \medskip\noindent
This means that the function\/ $\rho^{}_F(f\mu,f\nu)$ is additive with respect to $f$.
\end{theorem}

Generally, a real-valued measure $\nu$ is decomposed into three components
\begin{equation}
\label{e-decomp-nu}
 \nu =\nu_a +\nu^+_s +\nu^-_s,
\end{equation}
where $\nu_a$ is absolutely continuous, $\nu^+_s$ is positive and singular, and $\nu^-_s$ is negative and singular (with respect to $\mu$).

The  next result describes  the corresponding decomposition of $\rho^{}_F(\mu,\nu)$.

\begin{theorem} \label{..theorem_1a}
Let\/ $\mu$, $\nu$ be finite measures on a measurable space\/ $(X,\mathfrak A)$, where\/ $\mu$ is nonnegative and\/ $\nu$ is real-valued. Then
\begin{equation} \label{,,theorem_1a_sum}
  \rho^{}_F(\mu,\nu) = \rho^{}_F(\mu,\nu_a) +\rho^{}_F(0,\nu^+_s) +\rho^{}_F(0,\nu^-_s),
\end{equation}
where each term may be finite or equal to\/ $+\infty$, and
\begin{align} \label{,,theorem_1a_s+}
 \rho^{}_F(0,\nu^+_s) \pin&=\pin \nu^+_s(X)F'(+\infty),  \\[3pt] \label{,,theorem_1a_s-}
 \rho^{}_F(0,\nu^-_s) \pin&=\pin \nu^-_s(X)F'(-\infty).
\end{align}
Here we assume that if\/ $\nu^+_s =0$ or\/ $\nu^-_s =0$ then the corresponding product in the right-hand part of\/
%\eqref{,,theorem_1_sum}%
\eqref{,,theorem_1a_s+} or\/ \eqref{,,theorem_1a_s-} is zero regardless of the $($may be infinite\/$)$ multiplier\/ $F'(\pm\infty)$.
% ,  \\[3pt] \label{,,theorem_1_a}
% \rho^{}_F(\mu,\nu_a) \,&=\pin \int_X F\bigg(\frac{d\nu_a}{d\mu}\bigg)\pin d\mu.
\end{theorem}
\smallskip

There is quite a number of objects in analysis where one has to exploit not measurable partitions of unity but continuous ones (one of them will be considered in Section~\ref{s-n-t-entr}). To discuss this setting in our context we need the next definition.

%\begin{definition} \label{d-,,rho*}
Let $X$ be a topological space and $\mu,\, \nu$ be finite Borel measures, where $\mu$ is nonnegative and $\nu$ is real-valued. For any convex function $F\!:\mathbb R \to (-\infty,+\infty]$ set
\begin{equation} \label{,,rho*}
 \rho^{}_{F,c}(\mu,\nu) :=\, \sup_G\sum_{g\in G} \mu[g]\pin F\bigg(\frac{\nu[g]}{\mu[g]}\bigg),
\end{equation}
where supremum is taken over the set of all (finite) continuous partitions of unity\/ $G$ and we assume that if\/ $\mu[g]=0$ then the corresponding summand in the right-hand part is defined according to convention \eqref{,,72}.

\begin{theorem} \label{..theorem_3}
Let\/ $\mu$, $\nu$ be finite Borel measures on a metric space\/ $X$, where\/ $\mu$ is nonnegative and\/ $\nu$ is real-valued. Then for any convex lower semicontinuous function\/ $F$,
\begin{equation} \label{e,,theorem_3}
 \rho^{}_{F,c}(\mu,\nu) = \rho^{}_{F}(\mu,\nu).
\end{equation}
\end{theorem}

\vspace{-3pt}

\begin{remark} \label{..remark_1}
In fact instead of metrizability of $X$ in Theorem~\ref{..theorem_3} it suffices to require the density of the set of continuous functions $C(X)$ in the spaces $L^1(X,\mu)$ and $L^1(X,\nu)$ (which is always true for metrizable space $X$ or, as a variant, for \emph{regular} measures $\mu$, $\nu$).
\end{remark}

Now let us prove the above formulated results.

\medskip

\textit{Proof of Lemma~\ref{..71}.} If $s>0$ then \eqref{,,74} follows immediately from \eqref{,,73}. Note that \eqref{,,73} and \eqref{,,71} imply inequalities
\begin{equation} \label{,,lemma_2}
 F'(+\infty) \ge A, \qquad F'(-\infty) \le A.
\end{equation}

 \medskip\noindent
In turn along with \eqref{,,72} they imply \eqref{,,74}, provided $s=0$ and $x\ne 0$. Finally, in case $s=0$ and $x=0$ both sides in \eqref{,,74} become zero. \qed

\medskip

\textit{Proof of Lemma~\ref{..72}.}
 If $s,t>0$ then by convexity of $F$,
\begin{equation*}
 sF\bigg(\frac{x}{s}\bigg) +\pin tF\bigg(\frac{y}{t}\bigg) =\pin (s+t)\bigg(\frac{s}{s+t}F\bigg(\frac{x}{s}\bigg) +\frac{t}{s+t}F\bigg(\frac{y}{t}\bigg)\!\bigg)
 \ge\pin (s+t)F\bigg(\frac{x+y}{s+t}\bigg).
\end{equation*}

Consider the case when $s>0$ and $t=0$.

If $y=0$ then \eqref{,,75} turns into the equality $sF(x/s) =sF(x/s)$.

Suppose now that $y>0$. If at least one summand in the right-hand part of \eqref{,,75} is infinite, then \eqref{,,75} holds true. If both summands $sF(x/s)$ and $0F(y/0) =yF'(+\infty)$ are finite then the function $F$ must be finite and continuous on the whole interval $(x/s, +\infty)$. Hence in \eqref{,,75} one can pass to a limit as $t\to +0$ and obtain the desired inequality
\begin{equation*}
 sF\bigg(\frac{x+y}{s}\bigg) \pin\le\pin sF\bigg(\frac{x}{s}\bigg) + 0F\bigg(\frac{y}{0}\bigg).
\end{equation*}
The case $y<0$ is treated similarly.

It remains to analyse the case $s,t=0$ and $x,y\ne 0$. If $x$ and $y$ have the same sign (say $x,y>0$) then \eqref{,,75} turns into equality:
\begin{equation*}
 (x+y)F'(+\infty) \pin=\pin xF'(+\infty) +yF'(+\infty).
\end{equation*}
Suppose $x,\,y$ have different sings (say $x<0$ and $y>0$). Recall that $F'(+\infty) \ge F'(-\infty)$ (see \eqref{,,lemma_2}). Therefore in any case
\begin{equation*}
 xF'(-\infty) +yF'(+\infty) \,\ge\,
 \begin{cases}
     (x+y)F'(-\infty), & \hbox{if}\ \, x+y<0, \\[2pt]
     (x+y)F'(+\infty), & \hbox{if}\ \, x+y>0, \\[2pt]
     \pin 0, & \hbox{if}\ \, x+y=0,
 \end{cases}
\end{equation*}
which means that
\begin{equation*}
 0F\bigg(\frac{x}{0}\bigg) +\pin 0F\bigg(\frac{y}{0}\bigg) \pin\ge\, 0F\bigg(\frac{x+y}{0}\bigg).
\end{equation*}

 \medskip\noindent
Thus Lemma \ref{..72} is proved in all cases. \qed

\medskip

\textit{Proof of Lemma \ref{..lemma_rho}.} Consider a countable partition of unity $G = \{g_1,g_2,\dots\}$. First we prove that in this case the sum in \eqref{,,rho} is well-defined, i.\,e., that the limit
\begin{equation} \label{,,302}
 \lim_{n\to\infty} \sum_{i=1}^n \mu[g_i]\pin F\bigg(\frac{\nu[g_i]}{\mu[g_i]}\bigg)
\end{equation}
does exist, being either finite or equal to $+\infty$.

Set $h_n =\sum_{i\ge n} g_i$. Then by Levi's monotone convergence theorem,
\begin{equation} \label{,,303}
 \lim_{n\to\infty} \mu[h_n] =0 \quad \text{and}\quad \lim_{n\to\infty} |\nu| [h_n] =0,
\end{equation}
where $|\nu|$ denotes the total variation of $\nu$. Lemma \ref{..71} implies that
\begin{equation} \label{,,304}
 \sum_{i=n}^m \mu[g_i]\pin F\bigg(\frac{\nu[g_i]}{\mu[g_i]}\bigg) \ge\pin
 \sum_{i=n}^m \big(A\nu[g_i] +B\mu[g_i]\big) \ge\pin
  -|A||\nu|[h_n] -|B|\pin \mu[h_n].
\end{equation}
It follows from \eqref{,,303} and \eqref{,,304} that for any $\varepsilon>0$ there exists $N$ such that for all $n>N$ and $m\ge n$,
\begin{equation} \label{,,305}
 \sum_{i=n}^m \mu[g_i]\pin F\bigg(\frac{\nu[g_i]}{\mu[g_i]}\bigg) >\pin -\varepsilon.
\end{equation}

\medskip

Now we have two possibilities: if for any $\varepsilon>0$ there exists $N$ such that for all $n>N$ and $m\ge n$,
\begin{equation} \label{,,306}
 \sum_{i=n}^m \mu[g_i]\pin F\bigg(\frac{\nu[g_i]}{\mu[g_i]}\bigg) <\pin \varepsilon,
\end{equation}

 \medskip\noindent
then limit \eqref{,,302} does exist (being finite when all the summands in \eqref{,,302} are finite and equal to $+\infty$ when there is at least one infinite summand); otherwise, if assumption \eqref{,,306} fails, using its negation and \eqref{,,305} one can easily show that the limit still exists and equals $+\infty$.

Now let us check equivalence of finite and countable partitions for use in \eqref{,,rho}.

Each finite partition of unity $G$ in \eqref{,,rho} may be transformed into a countable one by adding countably many zero elements, so transition from finite to countable partitions cannot decrease the value of $\rho^{}_F(\mu,\nu)$. Thus, it suffices to proof that it cannot increase as well.

Let $\rho^{}_F(\mu,\nu)$ be defined by \eqref{,,rho} using countable partitions $G$. Then for any $c<\rho^{}_F(\mu,\nu)$ there exists a countable partition of unity $G=\{g_1,g_2,\dots\}$ such that
\begin{equation} \label{,,307}
 \sum_{i=1}^\infty\mu[g_i]\pin F\bigg(\frac{\nu[g_i]}{\mu[g_i]}\bigg) \,>\, c.
\end{equation}
Set $h_n =\sum_{i\ge n} g_i$. Combining Lemma \ref{..71} and \eqref{,,303} we obtain
\begin{equation} \label{,,308}
 \liminf_{n\to\infty} \mu[h_n]\pin F\bigg(\frac{\nu[h_n]}{\mu[h_n]}\bigg) \ge\pin
 \liminf_{n\to\infty} \big(A\nu[h_n] +B\mu[h_n]\big) =\pin 0.
\end{equation}
Consider a finite partition of unity $G_n =\{g_1,\dots,g_{n-1},h_n\}$. Now \eqref{,,307} and \eqref{,,308} imply
\begin{equation*}
 \sup_{n} \Bigg( \sum_{i=1}^{n-1} \mu[g_i]\pin F\bigg(\frac{\nu[g_i]}{\mu[g_i]}\bigg) +\pin
 \mu[h_n]\pin F\bigg(\frac{\nu[h_n]}{\mu[h_n]}\bigg)\!\Bigg) \,\ge\,
 \sum_{i=1}^\infty \mu[g_i]\pin F\bigg(\frac{\nu[g_i]}{\mu[g_i]}\bigg) >\pin c.
\end{equation*}
Then \eqref{,,307} is valid for some $G_n$ instead of $G$, which along with arbitrariness of the constant $c<\rho^{}_F(\mu,\nu)$ implies the statement of Lemma \ref{..lemma_rho}. \qed

\medskip

\textit{Proof of Theorem \ref{..73}.} If $g$ is an element of a measurable partition of unity $G$ then by Lemma \ref{..72},
\begin{equation*}
 \big(\mu_1[g] +\mu_2[g]\big) F\bigg(\frac{\nu_1[g] +\nu_2[g]}{\mu_1[g] +\mu_2[g]}\bigg) \pin\le\pin  \mu_1[g]\pin F\bigg(\frac{\nu_1[g]}{\mu_1[g]}\bigg) +\pin\mu_2[g]\pin F\bigg(\frac{\nu_2[g]}{\mu_2[g]}\bigg).
\end{equation*}
Summing this up over $g\in G$ and passing to suprema gives \eqref{,,76}. \qed

\medskip

\textit{Proof of Theorem \ref{..additivity}.} From Theorem \ref{..73} it follows that
\begin{equation*}
 \rho^{}_F\big((f_1+f_2)\mu,(f_1+f_2)\nu\big) \pin\le\pin
 \rho^{}_F(f_1\mu,f_1\nu) + \rho^{}_F(f_2\mu,f_2\nu).
\end{equation*}
So it suffices to prove the inverse inequality.

By definition, for any $c_i<\rho^{}_F(f_i\mu,f_i\nu)$, \,$i=1,2$, there exist measurable partitions of unity $G_i$, \,$i=1,2$, such that
\begin{equation} \label{,,84}
 \sum_{g\in G_i} \mu[f_ig]\pin F\bigg(\frac{\nu[f_ig]}{\mu[f_ig]}\bigg) >\pin c_i, \qquad i=1,2.
\end{equation}
For each $g\in G_i$ define the function
\begin{equation*}
 h_g =
 \begin{cases} \displaystyle
   f_ig/(f_1+f_2), & \text{if}\ \, f_1+f_2>0, \\[3pt]   \displaystyle
   g/2,          & \text{if}\ \, f_1+f_2=0.
 \end{cases}
\end{equation*}
Evidently, the collection  $H =\{\pin h_g\mid g\in G_1\cup G_2\pin\}$ forms a measurable partition of unity. Note that for each $g\in G_i$ we have the equality $(f_1+f_2)h_g =f_ig$. Therefore,
\begin{equation} \label{,,85}
 \sum_{h_g\in H} \mu[(f_1+f_2)h_g]\pin F\bigg(\frac{\nu[(f_1+f_2)h_g]}{\mu[(f_1+f_2)h_g]}\bigg) \pin=\pin \sum_{i=1}^2 \sum_{g\in G_i} \mu[f_ig]\pin F\bigg(\frac{\nu[f_ig]}{\mu[f_ig]}\bigg).
\end{equation}
From \eqref{,,84}, \eqref{,,85} it follows that
\begin{equation*}
 \rho^{}_F\big((f_1+f_2)\mu,(f_1+f_2)\nu\big) \pin>\pin c_1+c_2
\end{equation*}
and, by arbitrariness of $c_i<\rho^{}_F(f_i\mu,f_i\nu)$,
\begin{equation*}
  \rho^{}_F\big((f_1+f_2)\mu,(f_1+f_2)\nu\big) \pin\ge\pin
  \rho^{}_F(f_1\mu,f_1\nu) + \rho^{}_F(f_2\mu,f_2\nu). \qed
\end{equation*}

\medskip

\textit{Proof of Theorem \ref{..theorem_1a}.} The space $X$ can be decomposed into three disjoint measurable parts, say $X=X_a\sqcup X^+_s \sqcup X^-_s$, such that the measures $\mu$ and $\nu_a$ are supported on $X_a$ while $\nu^+_s$, $\nu^-_s$ are  respectively supported on $X^+_s$, $X^-_s$.
Denote by $f_a$, $f^+_s$, $f^-_s$ characteristic functions of these disjoint parts. Then
\begin{equation*}
 f_a\mu =\mu, \quad f_a\nu =\nu_a, \qquad f^+_s\mu =0, \quad  f^+_s\nu =\nu^+_s,
 \qquad f^-_s\mu =0, \quad f^-_s\nu =\nu^-_s,
\end{equation*}
and hence \eqref{,,theorem_1a_sum} follows from Theorem \ref{..additivity}.

%\begin{alignat*}{3}
% f_a\nu &=\nu_a,& \qquad f^+_s\nu &=\nu^+_s,& \qquad f^-_s\nu &=\nu^-_s, \\[3pt]
% f_a\mu &=\mu,&  f^+_s\mu &=0,&  f^-_s\mu &=0,
%\end{alignat*}

Proofs of equalities \eqref{,,theorem_1a_s+} and \eqref{,,theorem_1a_s-} are similar. For example,
\begin{equation*}
 \rho^{}_F(0,\nu^+_s) \pin=\pin \sup_G \sum_{g\in G} 0\pin F\bigg(\frac{\nu^+_s[g]}{0}\bigg) =\pin \sup_G\sum_{g\in G} \nu^+_s[g]F'(+\infty) \pin=\pin \nu^+_s(X)F'(+\infty). \qed
\end{equation*}

\medskip

To prove Theorem \ref{..theorem_3} we need the next

\begin{lemma} \label{l-part-cont-meas}
Let\/ $\mu$ be a positive finite Borel measure on a topological space $X$ such that\/
$C(X)$ is dense in\/ $L^1(X,\mu)$. Then for any measurable partition of unity\/ $G = \{g_1,\dots,g_n\}$ on\/ $X$ and any\/ $\varepsilon >0$ there exists a continuous partition of unity\/ $H =\{h_1,\dots,h_n\}$ on\/ $X$ such that\/ $\|h_i-g_i\| <\varepsilon$ in\/ $L^1(X,\mu)$ for all\/ $i\in \overline{1,n}$.
\end{lemma}

\textit{Proof.} Choose a small $\delta >0$ and approximate each $g_i$ by a continuous function $f_i$ satisfying $\|f_i-g_i\| <\delta$ in the space $L^1(X,\mu)$. Without loss of generality we can assume that the functions $f_i$ are strictly positive (which can always be guaranteed by replacing each $f_i$ by $\min\{f_i,0\} +\gamma$ with a small $\gamma >0$). Now define a continuous partition of unity with elements
\begin{equation*}
 h_i :=\frac{f_i}{\sum_{j=1}^n f_j}, \qquad i=1,\dots,n.
\end{equation*}
Clearly,
\begin{equation*}
 |h_i -f_i| \pin=\pin \bigg|\pin \frac{1-\sum_{j=1}^n f_j}{\sum_{j=1}^n f_j}\pin f_i\pin\bigg| \pin\leq\pin \bigg|\pin 1-\sum_{j=1}^n f_j\pin\bigg| \pin\leq\pin \sum_{j=1}^n |g_j -f_j|,
\end{equation*}

 \smallskip\noindent
which implies the estimate
\begin{equation*}
 \|h_i -g_i\| \leq \|h_i -f_i\| +\|f_i -g_i\| \leq n\delta +\delta
\end{equation*}
and by arbitrariness of $\delta$ finishes the proof Lemma~\ref{l-part-cont-meas}. \qed

 \medskip

\textit{Proof of Theorem \ref{..theorem_3}.} Since any continuous partition of unity is measurable it follows that $\rho^{}_{F,c}(\mu,\nu)\leq \rho^{}_F(\mu,\nu)$ and it is enough to prove the opposite inequality.

As in the proof of Theorem \ref{..theorem_1a} the space $X$ can be decomposed into three disjoint parts, $X=X_a\sqcup X^+_s \sqcup X^-_s$, such that the measures $\mu$ and $\nu_a$ are supported on $X_a$ while $\nu^+_s$, $\nu^-_s$ are  respectively supported on $X^+_s$, $X^-_s$. Denote by $f_a$, $f^+_s$, $f^-_s$ characteristic functions of these disjoint parts.

Theorem \ref{..theorem_1a} gives the following representation of $\rho^{}_F(\mu,\nu)$:
\begin{equation} \label{e-th-5:1}
 \rho^{}_F(\mu,\nu) \pin=\pin
 \sup_G\sum_{g\in G} \mu[g]\pin F\bigg(\frac{\nu_a[g]}{\mu[g]}\bigg) +\pin
 0F\bigg(\frac{\nu^+_s[f^+_s]}{0}\bigg) +\pin 0F\bigg(\frac{\nu^-_s[f^-_s]}{0}\bigg).
\end{equation}
Suppose for definiteness that $\nu^+_s[f^+_s] >0$ and $\nu^-_s[f^-_s] <0$ (otherwise the corresponding summands in the right-hand side of \eqref{e-th-5:1} may be omitted).

Note that in \eqref{e-th-5:1} one can assume that $\mu[g] >0$ for all $g$ since on the one hand the summands with $\mu[g] =0$ are equal to $0$ according to definition and on the other hand once
$\mu[g']=0$ and $\mu[g'']>0$ the pair $g'$, $g''$ can be replaced by one element $g = g' + g''$ in the partition $G$ that does not change the sum in \eqref{e-th-5:1} due to absolute continuity of $\nu_a$ with respect to $\mu$.

Now recalling lower semicontinuity of $F$ and definition of $F'(\pm\infty)$ the proof of theorem completes by applying Lemma~\ref{l-part-cont-meas} to partitions of unity $G' =\{\pin f_ag\mid g\in G\pin\}\cup\linebreak[2] \{f^+_s,f^-_s\}$ in the space $L^1(X,\mu+|\nu|)$. \qed \quad
%\textbf{[too hard paragraph?]}

\section{Sup\kern0.03em -sums principle for integral \emph{F}-divergence}  \label{s-1}

Here we present one of the principal results of the article uncovering interrelation between sup-sums $F$-divergence and integral $F$-divergence.

\begin{theorem} \label{..theorem_1} \hskip-0.15em \textup{(sup-sums principle)\,}
Let\/ $\mu$ and\/ $\nu$ be two finite measures on a measurable space\/ $(X,\mathfrak A)$, where\/ $\mu$ is non\-negative and\/ $\nu$ is real-valued, and\/ $\nu =\nu_a +\nu^+_s +\nu^-_s$ be decomposition\/~\eqref{e-decomp-nu}. Then
\begin{equation} \label{,,theorem_1_a}
 \rho^{}_F(\mu,\nu_a) \,=\pin \int_X F\bigg(\frac{d\nu_a}{d\mu}\bigg)\pin d\mu,
\end{equation}
and
\begin{equation} \label{,,theorem_1_sum}
  \rho^{}_F(\mu,\nu) = \pin \int_X F\bigg(\frac{d\nu_a}{d\mu}\bigg)\pin d\mu +  \nu^+_s(X)F'(+\infty) +  \nu^-_s(X)F'(-\infty).
\end{equation}

 \medskip\noindent
Here\/ $d\nu_a/d\mu$ denotes the Radon---Nikodym derivative and we assume that if\/ $\nu^+_s =0$ or\/ $\nu^-_s =0$ then the corresponding product in the right-hand part of\/ \eqref{,,theorem_1_sum} is zero regardless of the $($may be infinite\/$)$ multiplier\/ $F'(\pm\infty)$.
\end{theorem}

\begin{corollary} \label{..corollary}
For any\/ $f\in L^1(X,\mu)$ and any convex function\/ $F\!:\mathbb R \to (-\infty,+\infty],$
\begin{equation} \label{,,corollary}
 \int_X F(f)\,d\mu \,=\, \sup_G \sum_{g\in G} \mu[g]\pin F\bigg(\frac{\mu[fg]}{\mu[g]}\bigg),
\end{equation}
where supremum is taken over all measurable partitions of unity\/ $G$.
\end{corollary}

\textit{Proof of Theorem \ref{..theorem_1}.} Note that \eqref{,,theorem_1_sum} follows from \eqref{,,theorem_1_a} along with Theorem~\ref{..theorem_1a}.

Let us check that for each (no matter finite or countable) measurable partition of unity~$G$,
\begin{equation} \label{,,310}
 \sum_{g\in G} \mu[g]\pin F\bigg(\frac{\nu_a[g]}{\mu[g]}\bigg) \pin\le\pin \int_X F\bigg(\frac{d\nu_a}{d\mu}\bigg)\pin d\mu
\end{equation}

 \medskip\noindent
holds true. Indeed,
\begin{gather} \notag
 \sum_{g\in G} \mu[g]\pin F\bigg(\frac{\nu_a[g]}{\mu[g]}\bigg) \pin=\pin
 \sum_{\mu[g]>0} \mu[g]\pin F\bigg(\!\pin\int_X \frac{g}{\mu[g]}\, d\nu_a\bigg) \pin=\pin
 \sum_{\mu[g]>0} \mu[g]\pin F\bigg(\!\pin\int_X \frac{g}{\mu[g]}\pin \frac{d\nu_a}{d\mu}\, d\mu\bigg) \\[6pt]
 \le\pin \sum_{\mu[g]>0}\mu[g]\int_X \frac{g}{\mu[g]}\pin F\bigg(\frac{d\nu_a}{d\mu}\bigg)\pin d\mu
 \,= \pin \int_X \sum_{g\in G} g\pin F\bigg(\frac{d\nu_a}{d\mu}\bigg)\pin d\mu \,=\pin
 \int_X F\bigg(\frac{d\nu_a}{d\mu}\bigg)\pin d\mu, \label{,,simple_direction}
\end{gather}
where we exploited Jensen's inequality for the probability measure $(g/\mu[g])\pin d\mu$ and recall also that by convention \eqref{,,72} and absolute continuity of $\nu_a$ all the summands with $\mu[g] =0$ are zero.

From \eqref{,,310} it follows that the left-hand part in \eqref{,,theorem_1_a} does not exceed the right-hand one, and to finish the proof of Theorem \ref{..theorem_1} we have to verify the inequality
\begin{equation} \label{,,inverse}
 \sup_G\sum_{g\in G} \mu[g]\pin F\bigg(\frac{\nu_a[g]}{\mu[g]}\bigg) \pin\ge\pin
 \int_X F\bigg(\frac{d\nu_a}{d\mu}\bigg)\pin d\mu.
\end{equation}

For the convex function $F$ under consideration there exists a partition of real axis by three points $-\infty\le a\le b\le c\le +\infty$ (where everywhere there could be equalities) such that

\smallskip

i) \,$F(y) =+\infty$ for $y<a$ and $y>c$;

\smallskip

ii) \,$F(y)$ is nonincreasing, finite and continuous on $(a,b)$;

\smallskip

iii) \,$F(y)$ is nondecreasing, finite and continuous on $(b,c)$.

\smallskip

Let us decompose $X$ into seven subsets
\begin{equation} \label{,,311}
 X_{<a}, \quad X_a, \quad X_{(a,b)}, \quad X_b, \quad X_{(b,c)}, \quad X_c,\quad X_{>c},
\end{equation}
defined, respectively, by the conditions
\begin{equation*}
 \frac{d\nu_a}{d\mu}(x) < a, \quad \frac{d\nu_a}{d\mu}(x) =a, \quad a< \frac{d\nu_a}{d\mu}(x) <b,
 \quad \dots, \quad \frac{d\nu_a}{d\mu}(x) =c, \quad \frac{d\nu_a}{d\mu}(x) >c.
\end{equation*}
Some of them (and even the majority of them) may be empty; for example, if the function $F$ decreases on $(a,c)$, then $b=c$ and $X_{(b,c)} =\emptyset$, and if $F$ is finite everywhere then the sets $X_{<a}$, $X_a$, $X_c$, $X_{>c}$ will be empty.

Evidently, it is enough to prove inequality \eqref{,,inverse} for each of the sets \eqref{,,311} separately and then sum the components. In doing so partitions of unity $G$ on these sets should also be defined separately.

For the sets $X_{<a}$, $X_a$, $X_b$, $X_c$, $X_{>c}$ (some of them may by empty) inequality \eqref{,,inverse} is verified easily: it is sufficient to take a trivial partition $G$ consisting of a single unit function on the set considered.

Now consider the set $X_{(a,b)}$. Let us take an arbitrary number $\varepsilon>0$ and set
\begin{gather} \label{,,313}
 Y_i =\bigl\{ y\in (a,b)\bigm| i\varepsilon\le F(y) <i\varepsilon+\varepsilon\bigr\}, \qquad
 i\in\mathbb Z,\\[3pt]  \label{,,314}
 X_i = \Big\{ x\in X \Bigm| \frac{d\nu_a}{d\mu}(x) \in Y_i \pin\Big\}, \qquad i\in\mathbb Z.
\end{gather}
Clearly the sets $X_i$ form a partition of $X_{(a,b)}$ and their characteristic functions (that we denote by $g_i$) form a measurable partition of unity on $X_{(a,b)}$.

Note that by monotonicity of $F$ on $(a,b)$ the sets $Y_i$ are convex. Therefore, if $\mu(X_i)>0$, then
\begin{equation*}
 \frac{\nu_a(X_i)}{\mu(X_i)} \,=\, \frac{1}{\mu(X_i)}\int_{X_i} \frac{d\nu_a}{d\mu}\,d\mu \,\in\, Y_i,
\end{equation*}

 \medskip\noindent
and by definition of $Y_i$ we have
\begin{equation} \label{,,315}
 i\varepsilon \,\le\, F\bigg(\frac{\nu_a(X_i)}{\mu(X_i)}\bigg) <\, i\varepsilon +\varepsilon.
\end{equation}

Now \eqref{,,313}, \eqref{,,314}, \eqref{,,315} imply that
\begin{gather*}
 \int_{X_{(a,b)}} F\bigg(\frac{d\nu_a}{d\mu}\bigg)\pin d\mu \,=\pin
 \sum_{i\in\mathbb Z} \int_{X_i} F\bigg(\frac{d\nu_a}{d\mu}\bigg)\pin d\mu \,\le\pin
 \sum_{i\in\mathbb Z} (i\varepsilon+\varepsilon)\pin\mu(X_i) \\[6pt]
 \le\pin \sum_{i\in\mathbb Z} \mu(X_i) \bigg(F\bigg(\frac{\nu_a(X_i)}{\mu(X_i)}\bigg)
 +\,\varepsilon\bigg) \pin=\pin \sum_{i\in\mathbb Z} \mu[g_i] \pin
 F\bigg(\frac{\nu_a[g_i]}{\mu[g_i]}\bigg) +\, \varepsilon\mu\big(X_{(a,b)}\big).
\end{gather*}
By arbitrariness of $\varepsilon$ this implies inequality \eqref{,,inverse} for the set $X_{(a,b)}$.

For the set $X_{(b,c)}$ it is verified in the same way. Thus, Theorem \ref{..theorem_1} is proved. \qed

\medskip

\textit{Proof of Corollary \ref{..corollary}.} Take $\nu$ such that $d\nu/d\mu =f$ in \eqref{,,theorem_1_a}. \qed

\section{Sup\kern0.03em -sums principle for Kullback---Leibler\\
          divergence etc} \label{s-K-L}

If $\mu$ and $\nu$ are probability measures on $(X,\mathfrak A)$ and $\mu$ is absolutely continuous with respect to $\nu$ then Kullback---Leibler divergence $D_{\mathit{KL}}$ is defined as \begin{equation} \label{e.KL}
 D_{\mathit{KL}} (\mu\Vert \nu) := \int_X \ln\biggl(\frac{d\mu}{d\nu}\biggr)\, d\mu .
\end{equation}

The results of the foregoing section lead to the next

\begin{theorem} \label{t.KL}
Under the above conditions on\/ $\mu$ and\/ $\nu$,
\begin{gather} \label{e.KL0}
 D_{\mathit{KL}} (\mu\Vert \nu) \pin=\pin D_{\mathit{KL}} (\mu\Vert \nu_a) \pin = \pin \rho_{-\ln} (\mu, \nu) \pin = \pin \rho_{-\ln} (\mu, \nu_a) \\[9pt]  \label{e.KL0*}
 =\, \sup_G\sum_{g\in G} \mu[g] \ln\biggl(\frac{\mu[g]}{\nu[g]}\biggr) \pin =\,\sup_G\sum_{g\in G} \mu[g] \ln\biggl(\frac{\mu[g]}{\nu_a[g]}\biggr),
\end{gather}
where\/ $\nu_a$ is the absolutely continuous component of\/ $\nu$ with respect to\/ $\mu$ and supremum is taken over all\/ $($finite or countable\/$)$ measurable partitions of unity\/ $G$ on\/ $X$ and we assume that if\/ $\mu[g]=0$ then the corresponding summand in the sums  vanishes regardless of the second multiplier\/ $\ln(\mu[g]/\nu[g])$ or\/ $\ln(\mu[g]/\nu_a[g])$.
\end{theorem}

\textit{Proof.} According to \eqref{,,71}, we have $-\ln'(+\infty) =0$. Hence by Theorem \ref{..theorem_1},
\begin{equation} \label{,,kl_1}
 \rho_{-\ln}(\mu,\nu) \pin=\pin \rho_{-\ln}(\mu,\nu_a) \pin=
 \int_{X} -\ln\biggl(\frac{d\nu_a}{d\mu}\biggr)\,d\mu.
\end{equation}
It is easily seen that outside a set of zero measure $\mu$,
\begin{equation*}
 0 \pin<\pin \frac{d\mu}{d\nu_a} \pin=\pin \frac{d\mu}{d\nu} \pin<\pin +\infty.
\end{equation*}
Therefore,
\begin{equation} \label{,,kl_2}
 \int_X \ln\biggl(\frac{d\mu}{d\nu}\biggr)\,d\mu \pin=
 \int_X \ln\biggl(\frac{d\mu}{d\nu_a}\biggr)\,d\mu \pin=
 \int_{X} -\ln\biggl(\frac{d\nu_a}{d\mu}\biggr)\,d\mu.
\end{equation}

 \medskip\noindent
From \eqref{,,kl_1}, \eqref{,,kl_2} we obtain equalities \eqref{e.KL0}.

Recall that $\mu$ is absolutely continuous with respect to $\nu$ and hence with respect to $\nu_a$ as well. So if $\mu[g]\ne 0$ then $\nu[g]\ne 0$ and $\nu_a[g]\ne 0$ for any element $g$ of a measurable partition of unity on $X$. From this and definition \eqref{,,rho} of $\rho_{-\ln}(\mu,\nu)$ it follows that
\begin{equation*}
 \rho_{-\ln}(\mu,\nu) \pin=\,
 \sup_G\sum_{g\in G} \mu[g] \biggl(-\ln\biggl(\frac{\nu[g]}{\mu[g]}\biggr)\!\biggr) =\,
 \sup_G\sum_{g\in G} \mu[g] \ln\biggl(\frac{\mu[g]}{\nu[g]}\biggr),
\end{equation*}
where all summands with $\mu[g] =0$ are supposed to be zero. The analogous equality for $\rho_{-\ln}(\mu,\nu_a)$ may be got in the same way. Thus Theorem \ref{t.KL} is proved. \qed

\begin{remark} \label{r.KL}
The theorem just proved along with formula \eqref{,,kl_2} naturally suggests an extension of the definition of Kull\-back---Leibler divergence onto measures that are neither necessarily probability ones, nor mutually absolutely continuous. Namely, for any finite positive measures $\mu$, $\nu$ on a measurable space $(X,\mathfrak A)$ let us set
\begin{equation} \label{,,theorem_1-KL}
 D_{\mathit{KL}}(\mu\Vert \nu) := \int_X -\ln\biggl(\frac{d\nu_a}{d\mu}\biggr)\, d\mu.
\end{equation}
%Here $D_{\mathit{KL}} (\mu\Vert \nu)$ may be finite or equal to\/ $+\infty$.

The reasoning from the proof of Theorem~\ref{t.KL} shows that $D_{\mathit{KL}} (\mu\Vert \nu)$ defined in this way satisfies equalities \eqref{e.KL0} and \eqref{e.KL0*} as well.

%By the reasoning from the proof of Theorem~\ref{t.KL} for $D_{\mathit{KL}} (\mu\Vert \nu)$ defined in this way formulae \eqref{e.KL0} and \eqref{e.KL0*} are valid as well.
\end{remark}

\begin{remark} \label{r-cont-part}
If $X$ is a topological space and $\mu$ and $\nu$ are Borel measures such that the set $C(X)$ of continuous functions is dense in $L^1(X,\mu)$ and $L^1(X,\nu)$ (which is always true for a metrizable space $X$ or, as a variant, for \emph{regular} measures $\mu$, $\nu$) then recalling Theorem~\ref{..theorem_3} and Remark~\ref{..remark_1} one concludes that when applying \eqref{e.KL0} and \eqref{e.KL0*} to definition~\eqref{,,theorem_1-KL} we can equally use continuous (finite or countable) partitions of unity.
\end{remark}

\begin{remark}
As is known apart from Kullback---Leibler divergence many common divergences are special cases of $F$-divergence, corresponding to a suitable choice of $F$. For example, Hellinger distance corresponds to the function $F(t) =1-\sqrt{t}$, total variation distance corresponds to $F(t) =|t-1|$, Pearson $\chi^2$-divergence corresponds to $F(t) =(t-1)^2$, and for the function $F(t) =(t^\alpha -t)/(\alpha^2 -\alpha)$ we obtain the so-called $\alpha$-divergence.

Thus by choosing the corresponding convex functions $F$ one can write out the `sup-sums principles' of Theorem~\ref{t.KL} type for them.
\end{remark}

\begin{remark}
In the paper \cite{Sokol} the result of Theorem~\ref{t.KL} type was established for a sigma-finite measure $\nu$ and an absolutely continuous with respect to it measure $\mu$.
\end{remark}

\section{New definition for \emph{t}-entropy} \label{s-n-t-entr}

In this section we obtain on the base of Theorems~\ref{..theorem_3}, \ref{..theorem_1} and \ref{t.KL} a new transparent definition for $t$-entropy and clarify its relation to Kullback---Leibler divergence.

The $t$-entropy (its thorough definition see below) is a principal object of spectral analysis of operators associated with dynamical systems. In particular, in the series of articles \cite{ABL-2001,ABLS,ABL-2005,Bakh,ABL-Varp,ABL2} a relation between $t$-entropy and spectral radii of the corresponding operators has been established. Namely, it has been uncovered that $t$-entropy is the Fenchel---Legendre dual to the spectral exponent of operators in question.

For transparency of presentation let us recall the mentioned objects and results.

Hereafter $X$ is a Hausdorff compact space, $C(X)$ is the algebra of continuous functions on~$X$ taking real values and equipped with the max-norm, and $\alpha\!:X\to X$ is an arbitrary continuous mapping. The corresponding dynamical system will be denoted by $(X,\alpha)$.

Recall that a \emph{transfer operator} $A\!:C(X)\to C(X)$, associated with a given dynamical system is defined in the following way:

a) $A$ is a positive linear operator (i.\,e., it maps nonnegative functions to nonnegative ones); and

b) the following \emph{homological identity} for $A$ is valid:
\begin{equation} \label{1,,1}
 A(g \circ \alpha \cdot f) = gAf, \qquad g,f\in C(X).
\end{equation}
As an important and popular example of transfer operators one can take say the classical Perron---Frobenius
 operator, that is, the operator having the form
\begin{equation*}
 Af(x):= \sum_{y\in \alpha^{-1}(x)}a(y)f(y),
\end{equation*}
where $a\in C(X)$ is fixed. This operator is well defined when $\alpha$ is a local homeomorphism.

Let $A$ be a certain transfer operator in $C(X)$. In what follows we denote by $A_\varphi$ the family
of transfer operators in $C(X)$ given by the formula
\begin{equation*}
 A_\varphi f :=A(e^\varphi f), \qquad \varphi\in C(X).
\end{equation*}
Next, we denote by~$\lambda(\varphi)$ the \emph{spectral potential} of $A_\varphi$, namely,
\begin{equation*}
 \lambda(\varphi) := \lim_{n\to\infty}\frac{1}{n}\ln \Vert{A_\varphi^n}\Vert = \ln (r(A_\varphi)),
\end{equation*}
where $r(A_\varphi)$ is the spectral radius of operator $A_\varphi$.

We denote by $M(X)$ the set of all probability Borel measures on $X$. Recall that a measure $\mu\in M(X)$ is called \emph{$\alpha$-invariant} iff $\mu(g) =\mu(g\circ\alpha)$ for all $g\in C(X)$. The family of $\alpha$-invariant probability measures on $X$ is denoted by $M_\alpha(X)$.

The $t$-entropy is a certain functional on $M(X)$ denoted by $\tau(\mu)$ (its detailed definition will be given below).

The substantial importance of $t$-entropy is clearly demonstrated by the following variational
principle.

\begin{theorem} \label{1..1} \hskip -0.15em {\rm (\cite{ABL-Varp}, Theorem 5.6)\,}
Let\/ $A\!: C(X)\to C(X)$ be a transfer operator for a continuous mapping\/ $\alpha\!:X\to X$ of a
compact Hausdorff space\/~$X$. Then
\begin{equation*}
 \lambda(\varphi) \pin= \max_{\mu\in M_\alpha(X)} \bigl(\mu[\varphi] +\tau(\mu)\bigr), \qquad \varphi\in C(X).
\end{equation*}
\end{theorem}

One vividly notes the resemblance of this theorem to the Ruelle---Walters variational principle for the topological pressure \cite{Rue73,Wal75} uncovering its relation with Kolmogorov---Sinai entropy.

Among the principal ingredients in the proofs of the results leading to Theorem~\ref{1..1} is the so called `entropy statistic theorem'. This theorem plays in the spectral theory of weighted shift and transfer operators the role analogous to Shannon---McMillan---Breiman theorem in information theory \cite{McMill,AlCov} and its important corollary known as `asymptotic equipartition property' \cite[p.~135]{Bill}. The variational principles containing $t$-entropy and the objects therein serve as key ingredients of the thermodynamical formalism (see \cite{ABLS}, \cite{ABL2}, \cite{Lopes} and the sources quoted there).

Being so important $t$-entropy at the same time is rather sophisticated object to calculate.
The description of \mbox{$t$-entropy} not leaning on Fenchel---Legendre duality is not elementary and it took a substantial time and effort to obtain its `accessible' definition.

Namely, originally $t$-entropy $\tau(\mu)$ was defined in three steps (see, for example, \cite{ABL-Varp}).

\begin{definition} \label{d-t-entr1}
Firstly, for a given $\mu\in M(X)$, any $n\in\mathbb N$, and any continuous partition of unity $G =\{g_1,\dots,g_k\}$ we set
\begin{equation} \label{1,,3a}
 \tau_n(\mu,G) \pin:=\sup_{m\in M(X)}\sum_{g_i\in G}\mu[g_i] \ln\frac{m[A^ng_i]}{\mu[g_i]}.
\end{equation}
Here, if $\mu[g_i] = 0$ for some $g_i\in G$ then the corresponding summand in~\eqref{1,,3a} is
assumed to be zero regardless the value~$m[A^ng_i]$; if $m[A^ng_i] = 0$ for some $g_i\in G$ and at
the same time $\mu[g_i]>0$, then $\tau_n(\mu,G) = -\infty$.

Secondly, we put
\begin{equation} \label{1,,2a}
 \tau_n(\mu) :=\, \inf_G\tau_n(\mu,G);
\end{equation}

 \medskip\noindent
here the infimum is taken over all continuous partitions of unity $G$ in~$C(X)$.

And finally, the \emph{$t$-entropy} $\tau(\mu)$ is defined as
\begin{equation*} %\label{1,,2b}
 \tau(\mu) :=\pin \inf_{n\in\mathbb N}\frac{\tau_n(\mu)}{n}.
\end{equation*}
\end{definition}

Recently it was uncovered that for $\mu\in M_\alpha(X)$ (note that only such measures are essential for Theorem~\ref{1..1}) $t$-entropy could be defined in two steps \cite{BL2017}.

\begin{definition} \label{d-t-entr2}
First we set
\begin{equation} \label{e-tau-n}
 \tau_n(\mu):=\, \inf_G \sum_{g\in G} \mu [g]\ln\frac{\mu[A^n g]}{\mu[g]},
\end{equation}
where infimum is taken over the set of all continuous partitions of unity $G$ and we assume that if $\mu [g]=0$ then the corresponding summand in the right hand part of the equality is equal to $0$ independently of the value of $\mu[A^n g]$.

Now $\tau (\mu)$ is defined as
\begin{equation*} %\label{e-tau}
 \tau(\mu) :=\pin \inf_{n\in \mathbb N} \frac{\tau_n(\mu)}{n} .
\end{equation*}
\end{definition}

\medskip

In other words, in the original definition of $t$-entropy one should not calculate the supremum in
\eqref{1,,3a} but can simply put $m=\mu$ there. In \cite{BL2017} it was proved that this leads to the same value of $\tau_n(\mu)$ in \eqref{e-tau-n} as in \eqref{1,,2a}.

%to compare \eqref{1,,2a} and \eqref{e-tau-n}]

Of course, two steps are `better' (shorter) than three but even this two-steps definition of $t$-entropy is also rather sophisticated.

Note parenthetically that if one identifies a Borel measure $\mu$ on $X$ with a linear functional $\mu\!: C(X) \to \mathbb R$ given by
\begin{equation*}
 \mu [f] := \int_X f\, d\mu
\end{equation*}

 \medskip\noindent
then by Riesz's theorem there exists the only one regular Borel measure on $X$ defining the same functional. Thus, since in the foregoing definitions of $t$-entropy there were exploited only continuous functions (forming partitions of unity) we can assume that $t$-entropy is defined namely for regular measures $\mu$ (that are measures considered, in particular, in Theorem~\ref{..theorem_3} and Remark~\ref{..remark_1}).

The next theorem in essence gives a new definition of $t$-entropy and transparently establishes its relation to Kullback---Leibler divergence.

\begin{theorem} \label{t-t-entr-KL}
\hskip-0.15em \textup{($t$-entropy via Kullback---Leibler divergence)\,}
Let\/ $A$ be a transfer operator for a dynamical system\/ $(X, \alpha)$ then for any regular measure\/ $\mu \in M_\alpha(X)$ we have
\begin{equation*} %\label{e-tau-n-KL}
 \tau_n (\mu) \pin= \int_X \ln\frac{d(A^{*n}\mu)_a}{d\mu} \,d\mu \pin=\pin
 -D_{\mathit{KL}}(\mu\Vert A^{*n}\mu)
\end{equation*}
and
\begin{equation*} %\label{e-tau-KL}
 \tau(\mu) \pin=\pin
 \inf_{n\in \mathbb N}\,\frac{1}{n}\!\pin\int_X \ln\frac{d(A^{*n}\mu)_a}{d\mu}\,d\mu \pin=\pin -\sup_{n\in \mathbb N}\frac{D_{\mathit{KL}}(\mu\Vert A^{*n}\mu)}{n}\pin,
\end{equation*}

 \medskip\noindent
where\/ $A^*\!: C(X)^* \to C(X)^*$ is the operator adjoint to\/ $A$.
\end{theorem}

\textit{Proof.} Apply the reasoning of the proof of Theorem~\ref{t.KL} along with the reasoning of Remark~\ref{r-cont-part} to \eqref{e-tau-n}. Namely, set $\nu = A^{*n}\mu$ in this equality (so that $\mu[A^ng] = \nu [g]$) and  apply formulae \eqref{e.KL} -- \eqref{e.KL0*}. \qed

%%%%%%%%%%%%%%%%%%%%%%%%%%%%%%%%%%%%%%%%%%%%%%%%%%%%%%%%%%%%%%%%%%%%%%%%%%%%%%%%%%%%%%%%%%%%%%

\end{document}